\documentclass[graybox]{svmult}

\usepackage{type1cm}        % activate if the above 3 fonts are
\usepackage{makeidx}         % allows index generation
\usepackage{graphicx}        % standard LaTeX graphics tool
\usepackage{multicol}        % used for the two-column index
\usepackage[bottom]{footmisc}% places footnotes at page bottom

\usepackage{newtxtext}       % 
\usepackage[varvw]{newtxmath}       % selects Times Roman as basic font

\makeindex             % used for the subject index
\usepackage{mathtools}
\usepackage[colorlinks,citecolor=blue,urlcolor=blue]{hyperref}
\newcommand{\eps}{\varepsilon}
\newcommand{\Real}{\mathbb R}
\renewcommand{\E}{\mathsf{E}}
\renewcommand{\P}{\mathsf{P}}
\renewcommand{\Re}{\mathrm{Re}}
\renewcommand{\Im}{\mathrm{Im}}
\newcommand{\cov}{\mathsf{cov}}
\newcommand{\F}{\mathcal{F}}

\begin{document}

\motto{Dedicated to Prof. A.N.Shiryaev on the occasion of his 90-th birthday}

\title*{Asymptotic analysis in problems with fractional processes}
% Use \titlerunning{Short Title} for an abbreviated version of
% your contribution title if the original one is too long
\author{P. Chigansky  and M. Kleptsyna}
% Use \authorrunning{Short Title} for an abbreviated version of
% your contribution title if the original one is too long
\institute{P. Chigansky \at The Hebrew University, Mount Scopus, Jerusalem 91905, Israel \\ \email{Pavel.Chigansky@mail.huji.ac.il}
\and M. Kleptsyna \at Le Laboratoire Manceau de Math\'{e}matiques, Le Mans Universit\'{e}, France \\ \email{Marina.Kleptsyna@univ-lemans.fr}}
%
% Use the package "url.sty" to avoid
% problems with special characters
% used in your e-mail or web address
%
\maketitle

\abstract*{
Some problems in the theory and applications of stochastic processes can be reduced to solving integral equations. 
While explicit solutions for these equations are often elusive, valuable insights can be gained through their asymptotic analysis with respect to relevant parameters. This paper is a brief survey of some recent progress in the study of such equations related to processes with fractional covariance structure.
}

\abstract{
Some problems in the theory and applications of stochastic processes can be reduced to solving integral equations. 
While explicit solutions for these equations are often elusive, valuable insights can be gained through their asymptotic analysis with respect to relevant parameters. This paper is a brief survey of some recent progress in the study of such equations related to processes with fractional covariance structure.
}

\section{Introduction}

Linear integral equations frequently arise in the theory and applications of stochastic processes. 
In a typical setting, the quantity of interest is a particular functional of the solution to an integral equation, which involves the 
covariance operator of a stochastic process. Closed form solutions to such equations are rarely available.
One notable exception is a class of processes related to the Brownian motion. The covariance kernels for such processes can be  
identified with the Green function of an appropriate differential operator, easier to handle using the machinery of differential equations.
This class plays a significant role in physics and engineering since the Brownian motion is fundamental to modeling of the white noise in continuous time.

The purpose of this paper is to survey recent developments for processes that can be considered {\em fractional} analogues of the class mentioned above.
Just as the standard Brownian motion plays the defining role in that class, the main building block here is 
the fractional Brownian motion (fBm), that is, the centered Gaussian process $B^H=(B^H_t, t\in \Real_+)$ with the covariance function 
\begin{equation}\label{cov}
\E B^H_sB^H_t = \frac 1 2 \Big(|t|^{2H} +   |s|^{2H} - |t-s|^{2H} \Big),
\end{equation}
where $H\in (0,1)$ is its Hurst exponent. This process was introduced by A.Kol\-mo\-go\-rov in \cite{Kol40}, but its extensive study began only decades later, 
after publication of \cite{MvN68}. 
For $H=\frac 1 2$ the fBm reduces to the standard Brownian motion, but otherwise it is neither Markov, nor a semi-martingale.
It is the only self-similar Gaussian process with stationary increments  \cite[Theorem 1.3.3]{EM02}. For $H\in (\frac 1 2,1)$, the increments are positively 
correlated and have long range dependence:
$$
\sum_{n=2}^\infty \E  B^H_1  (B^H_n-B^H_{n-1}) =\infty.
$$
This property makes the fBm central in the study of natural phenomena exhibiting long memory \cite{PT17},  \cite{EM02}.

The theory of fBm is closely related to fractional analysis \cite{M08}, \cite{GMP23} and 
the method, reviewed in this paper, is applicable to a class of processes whose common feature is the fractional structure of  
covariance operator, similar to that of the fBm, see \cite[Section 4.4]{NP23} for typical examples. 
We will consider two standard types of problems associated with integral equations: the eigenvalue problem (Section \ref{sec:1})
and the equation of the second kind (Section \ref{sec:2}). Both are frequently encountered in problems of statistical inference,
optimal control and filtering of stochastic processes, see \cite{LS1, LS2}.
A typical result is the exact asymptotics of a relevant functional of the solution in a meaningful asymptotic regime. We will discuss a few concrete applications in which such results turn out to be useful. 
The main ideas behind the approach are sketched in Section \ref{sec:3}. 

\section{The eigenvalue problem}\label{sec:1}

The eigenvalue problem for the integral operator  
\begin{equation}\label{iop}
(K f)(t) = \int_0^T K(t,s)f(s)ds, \quad t\in [0,T],
\end{equation}
with a kernel function $K(\cdot,\cdot)$ consists of finding all pairs $(\lambda, \phi)$, the eigenvalues and eigenfunctions, which solve the 
equation 
\begin{equation}\label{eig}
K \phi = \lambda \phi.
\end{equation}
When the kernel is the covariance function $K(s,t) = \cov(X_s,X_t)$ of a stochastic process $X=(X_t, t\in [0,T])$, operator 
\eqref{iop} corresponds to the covariance operator of $X$. 
By the symmetry of $K(s,t)$, this operator is self-adjoint with respect to the usual inner product and, when, e.g., $K(\cdot,\cdot)\in L_2([0,T])$, 
it is also compact. Thus, by the Hilbert-Schmidt theorem, there are countably many eigenvalues $\lambda_n$, $n\in \mathbb N$, all of which are real and 
nonnegative with finite multiplicities. 
When put into the decreasing order, the sequence of  eigenvalues converges to zero as $n\to\infty$. The eigenfunctions $\phi_n$ corresponding to 
the eigenvalues, repeated according to multiplicity, can be chosen to form an orthonormal basis in $L_2([0,T])$. This basis is instrumental in many problems as it
diagonalizes the covariance operator.

An important implication is the Karhunen-Lo\`{e}ve series expansion 
of the process $X$ into $L_2$ convergent series of eigenfunctions, see \cite[Section, Ch 2, \S 2.3]{BTA04}, 
\begin{equation}\label{KL}
X_t = \sum_{n=1}^\infty \sqrt{\lambda_n}Z_n \phi_n(t)
\end{equation}
where  
$$
Z_n=\frac 1{\sqrt{\lambda_n}} \int_0^1 X(t)\phi_n(t)dt
$$  
are orthonormal random variables.
This expansion is useful for solving numerous problems  including characterization of equivalence  
of Gaussian measures \cite{Sh66}, small deviations problems \cite{NP23}, 
sampling from heavy tailed distributions \cite{VT13}, non-central limit theorems \cite{LRMT17}, nonparametric hypotheses testing \cite{N95}
and many others.

\begin{remark}
Both  eigenvalues $\lambda_n$ and  eigenfunctions $\phi_n(\cdot)$ depend on the interval length $T$. 
In this section, we will mainly discuss the results for the fBm, whose covariance function \eqref{cov}, due to its self-similarity, satisfies the scaling property 
$$
K(Tu,Tv)=T^{2H}K(u,v), \quad u,v\in[0,1].
$$
A simple change of variables in \eqref{eig} shows that in this case
\begin{equation}\label{lphi}
\lambda^{(T)}_n = T^{2H+1} \lambda^{(1)}_n, \quad \phi_n^{(T)}(t) = \frac 1 {\sqrt T} \phi_n^{(1)}(t/T), \quad t\in [0,T],
\end{equation}
where the superscript emphasizes the dependence on $T$. In view of this simple relation and for the sake of brevity, we will set $T=1$  
and omit it from the notations. All the formulas in this section can be readily adjusted to arbitrary $T$ using \eqref{lphi}.
For example, in the $L_2$-small ball probabilities asymptotics \eqref{L2sb}, for arbitrary $T>0$, the exponent $\beta(H)$ in \eqref{bg} 
should be multiplied by $T^{(2H+1)/2H}$ while $\gamma(H)$ remains intact.
\end{remark}

The eigenvalue problem \eqref{eig} can be solved explicitly only for a few processes. A prototypical example is the Brownian motion with $K(s,t)=s\wedge t$
for which there is a simple reduction to a boundary value problem for a linear ODE. An explicit solution to this ODE  gives the well-known exact formulas:
\begin{equation}\label{Bm}
\lambda_n = \frac 1 {\nu_n^2}\quad \text{and}\quad  \phi_n(t)=\sqrt 2 \sin (\nu_n t)\quad  \text{with} \quad \nu_n = (n-\tfrac 1 2)\pi, \quad n\in \mathbb N.
\end{equation}
More generally, reduction to boundary value problems for differential operators is possible for a whole class of the Green Gaussian processes \cite[Section 4]{NP23}, and in particular, for many processes related to the Brownian motion.  Consequently, the eigenvalue problem for such processes can be solved  exactly or, at least,  an accurate numeric approximation can be derived for its solutions. 

No such reduction is known for the covariance function \eqref{cov} of the fBm and it is unlikely that the eigenvalue problem can be solved explicitly 
in this case.
Numerical methods produce  reasonably accurate approximations for a few dozens of the first eigenvalues along with their eigenfunctions.
However, calculations tend to be unstable for smaller eigenvalues, for which asymptotic approximation is preferable, see \cite{VT13}, \cite{ChK} for some experiments.

The exact first order  asymptotics of  eigenvalues for the fBm was obtained in \cite{Br03a,Br03b}. It was shown that, for any $H\in (0,1)$ and arbitrarily small $\delta>0$,
\begin{equation}
\label{Bron}
\lambda_n =  \sin (\pi H)\Gamma(2H+1) (n\pi)^{-2H-1}  + o \big(n^{-r_H +\delta}\big) \quad   \text{as}\ n\to\infty,
\end{equation}
where $r_H : =  (2H+2)(4H+3)/(4H+5)$.
This asymptotic result was obtained by employing a specifically designed finite-dimensional approximation of the covariance operator. 
The leading order term in \eqref{Bron} was also derived by other means, see \cite{LP04} and \cite{NN04tpa}.

In some applications, such as those detailed below, it is crucial to know the exact second order term  in \eqref{Bron} and, moreover,   
to have an accurate asymptotic approximation for the eigenfunctions. It turns out that the eigenvalues and eigenfunctions in the fractional case have an asymptotic 
structure which, in a way, resembles \eqref{Bm}. It was shown in \cite{ChK} that
\begin{equation}\label{lambdan}
\lambda_n = \sin (\pi H)\Gamma(2H+1) \nu_n^{-2H-1},
\end{equation}
where the sequence $\nu_n$ satisfies 
\footnote{The second term on the right in \eqref{nun} is given in \cite{ChK} in a different, more cumbersome form}
\begin{equation}\label{nun}
\nu_n  = \big(n-\tfrac 1 2\big)\pi -   \frac{(H-\frac 12)^2 }{H+\frac 12}\frac \pi 2+ O(n^{-1}),\quad n\to \infty.
\end{equation}
This result improves the asymptotics in \eqref{Bron} up to the exact second order term and provides a tight estimate for the residual.  

The corresponding eigenfunctions are shown in \cite{ChK} to satisfy 
\begin{equation}\label{fBmeigfn}
\begin{aligned}
\phi_n& (t)  
= 
 \sqrt 2 \sin\big( \nu_{n}  t+\eta_H\big)  \\
&
 -    \int_0^{\infty}    
\Big(
e^{-  t \nu_n  u} f_0(u) +
(-1)^{n}   e^{-  (1-t) \nu_n  u}f_1(u)\
\Big)du + r_n(t) n^{-1} ,
\end{aligned}
\end{equation} 
where $\sup_n \|r_n\|_\infty<\infty$, the functions $f_0$ and $f_1$ have  explicit, albeit somewhat cumbersome, expressions and
$$
\eta_H = \frac 1 4 \frac{(H-\frac 1 2)(H-\frac 3 2)}{H+\frac 1 2}.
$$

Note that \eqref{fBmeigfn} is meaningful only if $\nu_n$ is known up to the $O(n^{-1})$ term, which is indeed the case in \eqref{nun}.
The second, integral term in \eqref{fBmeigfn} is responsible for the boundary layer: it vanishes asymptotically in the interior of 
the interval $[0,1]$ but persists at its endpoints, pushing the eigenfunctions to zero at $t=0$ and to certain specific values at $t=1$:
$$
\phi_n(1) = (-1)^{n} \sqrt{2H+1}   \big(1+O(n^{-1})\big), \quad n\to\infty.
$$
The boundary layer term may not be negligible for asymptotic approximation of various functionals of $\phi_n$, such as integrals, norms, etc.

As in the standard Brownian case, the approach taken in \cite{ChK} applies to other related process, such as the fractional Brownian noise \cite{ChK} 
(formal derivative of the fBm), the fractional Brownian bridge \cite{ChKM2} and its generalizations \cite{N19}, the mixed fractional Brownian motion \cite{ChKM20}, the fractional Ornstein-Uhlenbeck process \cite{ChKM-AiT}, the Riemann-Liouville process and some related fractional Sturm-Liouville problems \cite{ChK21}, the integrated fractional Brownian motion \cite{ChKM1}, \cite{ChKNR}. Typically, it produces an asymptotic approximation of the form similar to \eqref{lambdan} - \eqref{fBmeigfn} where all the ingredients are explicitly computable. The applicability to processes with variable Hurst parameter remains an interesting open question. Some results on the first order spectral asymptotics for such processes were recently obtained in \cite{KN21}, \cite{K22}. 

\medskip 

In the rest of this section we will discuss two problems to which such results are relevant: the $L_2$-small ball probabilities problem in Subsection \ref{ss:2.1} and the parameter estimation problem for the Ornstein-Uhlenbeck process in Subsection \ref{ss:2.2}.

\subsection{Small $\mathbf{L_2}$-small ball probabilities}\label{ss:2.1}
For a stochastic process $X=(X_t,t\in [0,1])$ and a norm $\|\cdot\|$, the small ball probability problem consists of finding the 
asymptotics of the probabilities 
$
\P(\|X\|\le \eps)
$ 
as $\eps\to 0$. This problem has a long history and numerous applications, see  \cite{LiS01}. In view of \eqref{KL}, the $L_2([0,1])$ norm satisfies the representation 
$$
\|X\|^2 = \sum_{n=1}^\infty \lambda_n Z_n^2. 
$$
When the process is Gaussian, the random variables $Z_n$ are standard normal,
which makes the problem more tractable. 

The exact asymptotic formula for the $L_2$ norm of the Brownian motion 
$B=(B_t, t\in [0,1])$ was derived in \cite{CM44}:
\begin{equation}
\label{CM}
\P(\|B\|_2\le \eps )=\frac{\ 4}{\sqrt \pi}\eps \exp\left(-\frac 1 8 \eps^{-2}\right)(1+o(1)), \quad \eps\to 0.
\end{equation}
The general solution in the Gaussian $L_2$ case was found in \cite{S74} in terms of the sequence of eigenvalues $\lambda_n$ of the covariance 
operator of $X$. 
This result was significantly simplified in \cite{DLL96} and \cite{NN04ptrf}, where the role of the asymptotic behavior of $\lambda_n$ was made more transparent. 
Essentially, it was shown that the leading order term in the asymptotic expansion of $\lambda_n$ determines the asymptotics of $\log \P(\|X\|_2\le \eps)$.  
If, in addition,  the second order term is known, the asymptotics of $\P(\|X\|_2\le \eps)$ can be found up to a multiplicative {\em distortion} constant.
A comprehensive account on $L_2$-small ball probabilities problem for Gaussian processes can be found in the recent survey \cite{NP23}.

Plugging the approximation \eqref{lambdan}-\eqref{nun} into the general formulas from \cite[Theorem 6.2]{NN04ptrf} extends \eqref{CM} to all values of $H\in (0,1)$: 
\begin{equation}\label{L2sb}
\P(\|B^H\|_2 \le \eps) = C(H) \eps^{\gamma(H)} \exp \left(-\beta(H)\eps^{-\frac 1 H}\right) (1+o(1))
, \quad \eps \to 0,
\end{equation}
where $C(H)$ is a constant which still remains unknown and  
\begin{equation}\label{bg}
 \beta(H)=  \frac{H}{(2H+1)^{\frac {2H+1}{2H}}}
\left(
\frac{\sin (\pi H) \Gamma(2H+1)}{\left(\sin  \frac{\pi}{2H+1} \right)^{2H+1}}
\right)^{\frac 1 {2H}}     \text{and}\ \  \gamma(H) =  \frac{(H-\frac 12)^2+1}{2H}.   
\end{equation}

\subsection{Drift estimation of the OU process}\label{ss:2.2}
A basic problem in statistical inference of stochastic processes is estimation of the drift coefficient of the 
Ornstein-Uhlenbeck (OU) process generated by the linear stochastic differential equation (SDE) 
\begin{equation}\label{OU}
dX_t = -\theta  X_t dt + dB_t,
\end{equation}
where $B =(B_t, t\in \Real_+)$ is the standard Brownian motion.  The true value of the parameter variable $\theta$, denoted by $\theta_0$, is unknown
and to be estimated from the sample $X^T=(X_t, t\in [0,T])$ up to a time $T>0$.
This problem has been the subject of extensive research over the years, see e.g. \cite{Nov72}, \cite[Chapter 17]{LS2}, \cite{Ku04}.

Denote by $\P^T_\theta$ the probability measure induced by the OU process on the space of continuous functions $C([0,T])$. 
By the Girsanov theorem, $\P^T_\theta$ is mutually absolutely continuous with respect $\P_0^T$, the Wiener measure, with the likelihood function 
\begin{equation}\label{like}
L_T(\theta) := \frac{d\P_\theta^T}{d\P_0^T}(X^T) = \exp \left(
-\theta\int_0^T X_tdX_t - \frac {1} 2 \theta^2 \int_0^T X_t^2dt 
\right ).
\end{equation}
Thus the maximum likelihood estimator (MLE) has the closed form formula 
$$
\widehat \theta_T = -\frac{\int_0^T X_tdX_t}{\int_0^T X_t^2dt}.
$$

The large sample asymptotics of $\widehat\theta_T$  depends on the sign of $\theta_0$. When $\theta_0>0$ the SDE \eqref{OU}
generates an ergodic process and, in this case, the MLE is consistent and asymptotically normal at rate $\sqrt T$:
\begin{equation}\label{OUAN}
\sqrt T (\widehat \theta_T - \theta_0) \xrightarrow[T\to\infty]{\mathcal {L}(\P_{\theta_0})} N\big(0, 2\theta_0\big).
\end{equation}
It can also be shown that $\widehat \theta_T$ is asymptotically efficient in the local minimax sense \cite[Ch 2, \S 11]{IKh}.
The asymptotics in the complementary non-ergodic case $\theta_0 \le 0$ is entirely different, see \cite[\S 3.5]{Ku04} for the details.

\subsubsection{Whittle's formula}
This result can be looked at from a different perspective. For $\theta>0$,  equation \eqref{OU} has the unique 
stationary solution with the spectral density 
\begin{equation}\label{sd}
f_\theta(\lambda) = \int_{-\infty}^\infty \E X_0X_t e^{-i\lambda t}dt = \frac 1{\lambda^2+\theta^2}, \quad \lambda \in \Real.
\end{equation}
The solution started at an arbitrary initial condition approaches the stationary one at an exponential rate. 
Hence, as it should be expected, the initial condition does not affect the estimation accuracy in the long run. 
The Fisher information in the stationary sample equals $TI(\theta)$, where the information rate $I(\theta)$ satisfies 
the continuous time analog of Whittle's formula
\cite{W53,W62}
\begin{equation}\label{W}
I(\theta) = \frac 1 {4\pi}\int_{-\infty}^\infty \left(\frac{\partial}{\partial \theta}\log f_\theta(\lambda)\right)^2d\lambda.
\end{equation}
For the particular spectral density in \eqref{sd} this gives $I(\theta)= (2\theta)^{-1}$, which 
makes asymptotics \eqref{OUAN} plausible: generically, the MLE should be asymptotically normal with the limit variance being equal to 
the reciprocal of the Fisher information rate, evaluated at the true value of the parameter.

This latter consideration suggests studying a more general OU type process, generated by  SDE \eqref{OU} with the Brownian motion replaced 
with a continuous Gaussian process $G$:  
\begin{equation}\label{OUG}
dX_t = -\theta  X_t dt + dG_t.
\end{equation}
If $G$ has stationary increments then, for $\theta>0$, we can  expect  that this equation has  unique stationary solution with the 
spectral density of the form, cf. \eqref{sd}, 
$$
f_\theta(\lambda) = \frac {g(\lambda)} {\lambda^2 +\theta^2},
$$
where $g(\lambda)$ is the ``spectral density'' of the ``derivative'' of $G$. If so, Whittle's formula suggests that asymptotics \eqref{OUAN}
should remain intact, regardless of the particular Gaussian process $G$. 

Making this heuristics exact is not a simple matter. The difficulty has to do with the subtlety of the absolute continuity relations 
between Gaussian measures on function spaces. In general, the likelihood function can be quite complicated and, consequently, 
neither the construction of MLE nor the validity of Whittle's formula \eqref{W} are a priori obvious, see \cite{DY83}.

\subsubsection{The innovation approach} If the process $X$ can be transformed 
into a semi-martingale by a causally invertible mapping, see \cite{HW77}, the problem becomes more tractable. 
Specifically, suppose there exists a kernel $g(s,t)$ such that the stochastic integral  
\begin{equation}\label{M}
M_t = \int_0^t g(s,t) dG_s, \quad t\ge 0, 
\end{equation}
defines a square integrable martingale with a strictly increasing (predictable) quadratic variation  
$\langle M\rangle_t$. If there also exists a kernel $h(s,t)$ so that 
$$
G_t = \int_0^t h(s,t)dM_s, \quad t\ge 0,
$$ 
then the natural filtrations of $M$ and $G$ coincide, $\F^M_t = \F^G_t$ for all $t\ge 0$.
In this case, $G$ is said to have the canonical innovation representation, see \cite[Ch. II]{HH76}, and the martingale 
$M$ is called fundamental for $G$.

By applying the same transformation to \eqref{OUG} we obtain a semi-martingale 
$$
Z_t = \int_0^t g(s,t) dX_s = -\theta \int_0^t X_s ds + M_t = -\theta \int_0^t Q(s) d\langle M\rangle_s + M_t,
$$
where we introduced an auxiliary process 
$$
Q(t) = \frac{d}{d\langle M\rangle_t}\int_0^t X_s ds.
$$
If the filtrations of $X$ and $Z$ coincide, $\F^X_t=\F^Z_t$ for all $t\ge 0$, then 
by Girsanov's theorem the likelihood function has the form, cf. \eqref{like},
$$
L_T(\theta) = \frac{d\P_\theta^T}{d\P_0^T}(X^T) = \exp \left(
-\theta\int_0^T Q(t) dZ_t - \frac {1} 2 \theta^2 \int_0^T Q(t)^2 d\langle M\rangle_t
\right ),
$$
and, consequently, the MLE is given by the formula
$$
\widehat \theta_T = -\frac{\int_0^T Q(t)dZ_t}{\int_0^T Q(t)^2 d\langle M\rangle_t}.
$$
The scaled estimation error satisfies 
$$
\sqrt T (\widehat \theta_T -\theta_0) = -\frac{\frac 1 {\sqrt T}\int_0^T Q(t)dM_t}{\frac 1 T\int_0^T Q(t)^2 d\langle M\rangle_t}.
$$
If we now manage  to prove that 
\begin{equation}\label{manage}
\frac 1 T\int_0^T Q(t)^2 d\langle M\rangle_t\xrightarrow[T\to\infty]{\P_{\theta_0}} \frac 1{2\theta_0},
\end{equation}
then, by the CLT for stochastic integrals \cite[Theorem 1.19]{Ku04},
$$
\frac 1 {\sqrt T}\int_0^T Q(t)dM_t\xrightarrow[T\to\infty]{\mathcal{L}(\P_{\theta_0})} N\left(0,\frac 1 {2\theta_0}\right),
$$
and  asymptotics \eqref{OUAN} follows. Applicability of this approach requires constructing a suitable canonical innovation representation. 

\subsubsection{Fractional Brownian motion}
The fBm is known to have such a representation,  see \cite{MG69}, \cite{NVV99} and \eqref{M} for $G:=B^H$ holds with 
\begin{equation}\label{gst} 
g(s,t) = \kappa_H^{-1} s^{\frac 1 2-H}(t-s)^{\frac 1 2 -H},
\end{equation}
where $\kappa_H$ is an explicit constant. The quadratic variation bracket of the fundamental martingale $M$ also has the simple closed form:
\begin{equation}\label{w}
\langle M\rangle_t = \lambda_H^{-1} t^{2-2H}
\end{equation}
with another explicit constant $\lambda_H$. Convergence \eqref{manage} was proved in \cite{KLeB02} by means of asymptotic analysis 
as $T\to\infty$ of the Laplace transform 
\begin{equation}\label{Laplace}
L_T(\mu) = \E_{\theta_0} \exp \left(-\mu \frac 1 T \int_0^T Q(s)^2 d\langle M\rangle_s \right), \quad \mu\in \Real.
\end{equation}
A closed form Cameron-Martin type formula, involving certain Riccati ODE, was derived for $L_T(\mu)$ in \cite{KLeB02}.
For the specific functions in \eqref{gst} and \eqref{w} its solutions can be expressed in terms of modified Bessel functions, 
which makes the asymptotic calculations feasible. 

\subsubsection{Mixed fractional Brownian motion}\label{sec:mfBm}
A large class of processes for which canonical innovation representation exists, consists of additive mixtures  
$$
G_t = B_t + V_t, \quad t\ge 0,
$$
of the Brownian motion $B$ and an independent centered Gaussian process $V$. 
It is known from \cite{Sh66} that the measure induced by $G^T=(G_t, t\in [0,T])$ is mutually absolutely continuous with respect to the Wiener measure if and only if its covariance satisfies 
$$
 \cov (V_s,V_t) = \int_0^t \int_0^s R(u,v)dudv
$$
with the kernel  
\begin{equation}\label{Shepp}
\int_0^T\int_0^T R(s,t)^2 dsdt <\infty.
\end{equation}
Under this condition a canonical innovation representation for $G$ was constructed in \cite{Hi68}.

A particular mixture of this type, the so-called mixed fBm,
\begin{equation}\label{mfBm}
G_t = B_t + B^H_t, \quad t\ge 0,
\end{equation}
arises in mathematical finance \cite{Ch01}, \cite{Ch03}, see also \cite{DSSV24}. Shepp's condition \eqref{Shepp} for the fBm covariance 
\eqref{cov} holds if and only if $H\in (\frac 34,1)$, in which case $G$ is a semi-martingale.
Thus, for such values of $H$, the canonical innovation representation is given in \cite{Hi68}, see also \cite{Ch03b}.

It turns out, however, that the mixed fBm has such a representation for all $H\in (0,1)$. As shown in \cite{CCK}  
the martingale $M_t = \E (B_t|\F^G_t)$ is fundamental and it satisfies \eqref{M}
with $g(s,t)$ being the unique continuous solution to the integro-differential equation 
\begin{equation}\label{geq}
g(s,t) + \frac \partial{\partial s} \int_0^t   H|s-r|^{2H-1}\mathrm{sign}(s-r) g(r,t)dr =1,\quad 0< s,t\le T, 
\end{equation}
where the kernel is the derivative of the fBm covariance \eqref{cov}. 
Moreover, cf. \eqref{w},
\begin{equation}\label{pqv}
\langle M \rangle_t = \int_0^t g(s,t) ds.
\end{equation}

Asymptotic analysis of the Laplace transform \eqref{Laplace} in  \cite{CKstat}  reveals that  \eqref{OUAN} holds if 
$\langle M\rangle_t$ satisfies the growth conditions:
\begin{equation}\label{A}
\frac 1 t \max \left(\frac{dt}{d\langle M\rangle_t}, \frac{d\langle M\rangle_t}{dt}\right)
\xrightarrow[]{t\to\infty} 0\quad \text{and}\quad
\int_0^\infty \left(\frac d {dt} \log \frac d {dt} \langle M\rangle_t\right)^2 dt <\infty.
\end{equation}
This is where the asymptotic approximation for the eigenvalues of the covariance operator comes into play.
Let us consider the simpler case $H>\frac 1 2$,
for which the derivative in \eqref{geq} can be taken under the integral. Thus a simpler, integral equation  is obtained:
\begin{equation}\label{geqred}
g(s,t) +  \int_0^t  c_H |s-r|^{2H-2} g(r,t) dr =1,\quad 0< s,t\le T, 
\end{equation}
with $c_H = H(2H-1)$, and, moreover, \eqref{pqv} simplifies to 
\begin{equation}\label{Mbr}
\langle M\rangle_t = \int_0^t g(t,t)^2dt.
\end{equation}
In order to verify conditions \eqref{A}, we need to estimate the growth of $g(t,t)$ as $t\to\infty$.

To this end, let us define the small parameter $\eps:= T^{1-2H}$ and rescale \eqref{geqred} accordingly, so that 
the function $u_\eps(x) := T^{2H-1}g(xT,T)$ solves the integral equation of the second kind 
\begin{equation}\label{2nd}
\eps u_\eps (x)  + \int_0^1 c_H |y-x|^{2H-2}u_\eps(y) dy = 1, \quad x\in [0,1],
\end{equation}
and the bracket in \eqref{Mbr} satisfies 
\begin{equation}\label{as}
\frac{d\langle M\rangle_T}{dT} = g^2(T,T) = \eps^2 u_\eps^2(1).
\end{equation}

Equation \eqref{2nd} has unique  solution, continuous on $[0,1]$. For $\eps =0$, it reduces to equation of the first kind  
\begin{equation}\label{1st}
 \int_0^1 c_H |y-x|^{2H-2} u_0(y)dy = 1, \quad x\in (0,1),
\end{equation}
which has the explicit solution, constructed in \cite{LC67}:
$$
u_0(x) = a_H x^{\frac 1 2 -H}(1-x)^{\frac 1 2-H}, \quad x\in (0,1),
$$
with some constant $a_H$. 
This function explodes at the endpoints of the interval and therefore $u_\eps(1)$ can be expected to diverge to infinity as $\eps\to 0$. 
In view of \eqref{as}, we need the exact rate at which this divergence occurs.

This rate can be found by expanding the solutions to \eqref{2nd} and \eqref{1st} into the Hilbert-Schmidt series 
$$
 u_\eps   =\sum_{n=1}^\infty \frac{\langle 1, \phi_n\rangle}{\eps+\lambda_n}\phi_n\quad \text{and}\quad 
 u_0  =\sum_{n=1}^\infty \frac{\langle 1, \phi_n\rangle}{ \lambda_n}\phi_n,
$$
where $\lambda_n$ and $\phi_n$ are the eigenvalues and eigenfunctions of the  kernel $K(x,y)=c_H|x-y|^{2H-2}$. This is the covariance 
the fractional Brownian "noise", the formal "derivative" of the fBm.
By subtracting \eqref{1st} from \eqref{2nd} we get 
\begin{equation}\label{expan}
\begin{aligned}
u_\eps(1) =\, & \eps^{-1} \int_0^1 c_H |y-1|^{2H-2} \big(u_0(y)-u_\eps(y)\big)dy = \\
& \sum_{n\ \text{odd}} 
\frac{\langle 1,\phi_n\rangle}{\lambda_n(\eps+\lambda_n)} \int_0^1 \phi_n(y) c_H(1-y)^{2H-2}dy,
\end{aligned}
\end{equation}
where for even indices $\langle 1,\phi_n\rangle=0$.
The growth of this series as $\eps\to 0$ can be estimated using the asymptotics of the eigenvalues from \cite{Ukai}:
$$
\lambda_n = c_1(H) n^{1-2H}(1+o(1)), \quad n\to \infty
$$
and of the scalar products from \cite{ChK}:
$$
\begin{aligned}
&  
\int_0^1 x^{-\beta} \phi_n(x)dx = c_2(H) n^{\beta-1} (1+o(1)),
\\
&
\int_0^1  \phi_{2n-1}(x)dx = c_3(H) n^{-\frac 1 2 -H} (1+o(1)),
\end{aligned}\qquad n\to\infty.
$$
where $c_j(H)$'s are nonzero constants and $\beta\in (0,1)$.
Substitution into  \eqref{expan} gives
$$
u_\eps(1) =  
C_1 \sum_{n=1}^\infty 
\frac{n^{-\frac 1 2 -H}}{ \eps+n^{1-2H}} (1+o(1))  = C_2 \eps^{-\frac 1 2} (1+o(1)), 
\quad \eps\to 0
$$
for some constants $C_1$ and $C_2$.
Consequently, $\eps^2 u_\eps^2(1) = C_2  \eps (1+o(1))$ as $\eps\to 0$ which translates to 
$$
\frac{d\langle M\rangle_T}{dT} =   C_2 T^{1-2H} (1+o(1)), \quad T\to \infty.
$$
This verifies the first condition in \eqref{A} and the second condition can be checked similarly, see \cite{CKstat}.

\section{Equations of the second kind}\label{sec:2} 

In this section we consider several problems which involve integral equations of the second kind  
\begin{equation}\label{kind2}
\eps g(t) + \int_0^T K(s,t)g(s)ds = f(t),\quad t\in [0,T],
\end{equation}
where $K(\cdot,\cdot)$ is the covariance function of a stochastic process,  $f(\cdot)$ is a given function
and $\eps$ and $T$ are real positive  parameters.
General theory of such equations establishes existence and uniqueness of solutions in  appropriate function spaces, \cite[Ch. IV]{RN55}. 
A typical objective is to compute a particular functional of the solution 
$g(\cdot)$, whose specific form depends on the problem.

In the context of asymptotic analysis, it is natural to study  such equations in two regimes. 
One is to fix $T>0$ and let $\eps\to 0$.  
Then \eqref{kind2} can be viewed as a small perturbation of the limiting equation of the first kind
$$
\int_0^T K(s,t)g^0(s)ds = f(t),\quad t\in [0,T],
$$
obtained by setting $\eps=0$. Let us denote the solution to 
 \eqref{kind2} by $g^\eps(\cdot)$, adding the superscript to emphasize the dependence on the relevant parameter. 
Such perturbation often turns out to be singular in the sense that $g^0(\cdot)$ does not 
provide a uniform approximation for $g^\eps(\cdot)$ on the whole interval as $\eps\to 0$. This is clearly demonstrated by the equations 
\eqref{2nd} and \eqref{1st}. Singularly perturbed integral equations have been studied for several types of
kernels, see, e.g.,  \cite{LS88, LS93}, \cite{Sh06} for rigorous theory and \cite{AO87b,AO87a} for more formal expansions. 

Another asymptotic regime is to fix $\eps>0$ and let $T\to \infty$. In this case, \eqref{kind2} can be viewed as the {\em finite section} approximation 
of  equation on the whole positive semiaxis 
$$
\eps g_\infty(t) + \int_0^\infty K(s,t)g_\infty(s)ds = f(t),\quad t\in \Real_+.
$$
For integrable difference kernels, i.e. when $K(s,t)=q(|s-t|)$ for some function $q\in L_1(\Real_+)$, explicit solution $g_\infty(\cdot)$ to such equations can be often constructed using the Wiener-Hopf method \cite{Krein58}.
The objective of the asymptotic analysis is to quantify its proximity to  $g_T(\cdot)$, the solution to \eqref{kind2}, see \cite{GF74}.

\medskip 
 
Below we discuss in detail two problems where such asymptotic regimes naturally arise: the filtering problem for the fractional counterpart of the Kalman-Bucy model in Subsection \ref{ss:KB} and the Hurst parameter estimation problem for the mixed fBm in Subsection \ref{ss:H}.

\subsection{The linear filtering problem}\label{ss:KB}

The filtering problem is concerned with estimation of signals from their noisy observations. 
A typical model in continuous time consists of a pair of random processes, the signal or state process $X=(X_t,\ t\in \Real_+)$ and 
the observation process $Y=(Y_t,\ t\in  \Real_+)$, and it is required to estimate the state $X_T$ at a time $T$ 
given the trajectory $Y^T =(Y_s, s\in [0,T])$. This  problem often arises in engineering applications and it 
has been studied extensively since World War II, see \cite{C14} for a concise historical account.

In a basic state-space model of R. Kalman and R. Bucy \cite{KB61} the state is the OU process 
generated by the linear SDE
\begin{equation}\label{Xeq}
X_t = \beta\int_0^t  X_s ds + W_t, \quad t\in \Real_+,
\end{equation}
where $\beta\in \Real$ is a constant and $W=(W_t,t\in \Real_+)$ is the Brownian motion. The state process is observed with the {\em additive Gaussian white noise}
so that $Y$ is given by 
\begin{equation}\label{Yeq}
Y_t = \mu\int_0^t  X_s ds +\sqrt{\eps} V_t, \quad t\in \Real_+,
\end{equation}
where $\mu\ne 0$ is another constant, $V=(V_t, t\in \Real_+)$ is an independent Brownian motion and $\eps>0$ is the observation noise intensity parameter.

The optimal estimator of $X_T$ is the conditional expectation $\widehat X_{T} = \E(X_T|\F^Y_T)$ given $\F^Y_T = \sigma\{Y_\tau, \tau \in [0,T]\}$.
It minimizes the mean squared error among all estimators based on the observed sample $Y^T$.
For the model \eqref{Xeq}-\eqref{Yeq} the optimal estimator $\widehat X_t$ satisfies the Kalman-Bucy filtering SDE 
$$
d\widehat X_t = \beta \widehat X_t dt + \frac{\mu P_t}{\eps} (dY_t - \mu \widehat X_t dt), 
$$
and the filtering error $P_t = \E(X_t-\widehat X_t)^2$ solves the Riccati ODE
\begin{equation}\label{Riccati}
\dot P_t = 2\beta P_t + 1 -  (\mu/\sqrt\eps)^2 P_t^2,
\end{equation}
subject to zero initial condition. 

A simple calculation shows that the 
filtering error has the steady-state limit 
\begin{equation}\label{PTinf}
\lim_{T\to\infty} P_T\Big(\beta, \frac \mu {\sqrt{\eps}}\Big) = \frac{\beta + \sqrt{\beta^2+\mu^2/\eps}}{\mu^2/\eps} 
\end{equation}
and scales with respect to the small noise intensity as
\begin{equation}\label{PTeps}
P_T\Big(\beta, \frac \mu {\sqrt{\eps}}\Big) = \frac{\sqrt{\eps}}{\mu}\big(1+o(1)\big),\quad \text{as\ } \eps  \to 0, \quad \forall T>0.
\end{equation} 
Such asymptotic quantities are of considerable interest, as they reveal the fundamental accuracy limitations
in the problem.

The Kalman-Bucy model can be considered in a greater generality, with $W$ and $V$ being independent Gaussian processes. 
Assuming that the processes are sufficiently regular, the optimal estimator in this case is given by the stochastic integral 
\begin{equation}\label{stint}
\widehat X_T = \frac 1{\mu} \int_0^T   g_T(s)dY_s,
\end{equation}
where $g_T(\cdot)$ solves the integro-differential equation 
\begin{equation}\label{maineq}
\frac {  \eps}{\mu^2}\frac{\partial}{\partial s} \int_0^T    g_T(r)\frac{\partial}{\partial r} K_V(r,s)dr
\\
+  \int_0^T K_X(r,s)   g_T(r)dr =  K_X(s,T),
\end{equation} 
with $K_V(s,t)$ and $K_X(s,t)$ being the covariance functions of the observation disturbance $V$ and the state process $X$.
The minimal  error can be computed by the formula   
\begin{equation}\label{PTfun}
P_T \Big(\beta, \frac{ \eps}{\mu^2}\Big)= \frac {\eps} {\mu^2} \left(
\frac{\partial}{\partial s} \int_0^T     g_T(r)\frac{\partial}{\partial r} K_V(r,s)dr
\right)_{\displaystyle\big|s:=T},
\end{equation}
where the dependence on the model parameters is made explicit in the notation. 

It is insightful to trace back how  general equations \eqref{maineq}-\eqref{PTfun} simplify in the Kalman-Bucy case.
When $V$ is the Brownian motion, i.e. $K_V(r,s)=r\wedge s$, the first term in \eqref{maineq} reduces to $g_T(s)$ and an integral 
equation of the second kind is obtained, cf. \eqref{kind2}. The functional in \eqref{PTfun} also takes the simpler, local form  
\begin{equation}\label{endpnt}
P_T \Big(\beta, \frac{ \eps}{\mu^2}\Big)= \frac {\eps} {\mu^2} g_T(T).
\end{equation}
Under mild regularity conditions, the solution $g_T(\cdot)$  belongs to $L_2([0,T])$ and consequently the stochastic integral in \eqref{stint} is well defined. Furthermore, when $W$ is the Brownian motion, the exponential structure of the covariance function $K_X(s,t)$ of the OU process 
allows to derive  the Riccati equation for the value of its solution at the endpoint, cf.  \eqref{endpnt} and \eqref{Riccati}.

Without such a structure, the general problem appears to be much more complex, see \cite{KLeB01}.
In view of relevance of the fBm to noise modeling,  \cite{BP88}, \cite{WG07}, \cite{PT17}, it makes sense to consider 
the fractional counterpart of the Kalman-Bucy model, when $W$ and $V$ in \eqref{Xeq} and \eqref{Yeq} are independent 
fBm's with Hurst parameters $H_1$ and $H_2$, respectively. 
Equation \eqref{maineq} in this case can hardly be expected to reduce to an ODE, yet alone to have a closed form solution.
Therefore asymptotic formulae in the spirit of \eqref{PTinf} and \eqref{PTeps} are of considerable interest as they can be used as 
benchmarks for design of simpler yet near optimal filtering procedures.

The first result in this direction was obtained in \cite{KlLB02}, where the steady-state error  was computed in the case of equal Hurst exponents $H:=H_1=H_2\in (\frac 1 2,1)$:
\begin{equation}\label{KLeB}
\lim_{T\to\infty}P_T \Big(\beta, \frac{ \eps}{\mu^2}\Big) = \frac{\Gamma(2H+1)}{2(\beta^2 + \mu^2/\eps)^H}\left(
1+ \sin (\pi H) \frac{\sqrt{\beta^2 + \mu^2/\eps}+\beta}{\sqrt{\beta^2 + \mu^2/\eps}-\beta}
\right).
\end{equation}
The small noise asymptotics can be readily derived from this formula:
\begin{equation}\label{Hsn}
P_\infty \Big(\beta, \frac{ \eps}{\mu^2}\Big) = \frac 1 2 \Gamma(2H+1)\big(1+\sin(\pi H)\big) (\eps/\mu^2)^H \big(1+o(1)\big), \quad \eps\to 0.
\end{equation}
For $H=\frac 1 2$, formula \eqref{KLeB} reduces to \eqref{PTinf} and asymptotics \eqref{Hsn} coincides with \eqref{PTeps}.

The limit \eqref{KLeB} is found in  \cite{KlLB02} by means of asymptotic analysis of the integro-differential Riccati equation from \cite{KLeB01}.
This approach does not seem to apply beyond the case of equal Hurst exponents $H_1= H_2$. It turns out, however, that equation \eqref{maineq} can be 
tackled directly  and, in fact, this reveals a more informative picture. As shown in \cite{AChKM22},  the steady state limit exists for all 
$H_1,H_2\in (0,1)$:
\begin{equation}\label{limerr}
P_\infty \Big(\beta, \frac{ \eps}{\mu^2}\Big) := \lim_{T\to\infty} P_T \Big(\beta, \frac{ \eps}{\mu^2}\Big),
\end{equation}
and it is related to the small noise asymptotics through the scaling, cf. \eqref{Hsn}:
\begin{equation}\label{epsas}
\lim_{\eps \to 0}\eps^{-\nu} P_T \Big(\beta, \frac{ \eps}{\mu^2}\Big) = P_\infty(0,\mu) \quad \text{where}\quad \nu = \frac{H_1}{1+H_1-H_2}.
\end{equation}

The limit error \eqref{limerr} has a rather interesting structure.
It is determined by a particular extension of the function  
\begin{equation}\label{spd}
f(\lambda) = \mu^2 \frac{\kappa(H_1)|\lambda|^{1-2H_1}}{\lambda^2+\beta^2}+\eps \kappa(H_2) |\lambda|^{1-2H_2}, \quad \lambda \in \Real,
\end{equation}
to the complex plane. In the upper half plane, it is defined  by the formula 
$$
\Lambda(z;H_1,H_2) :=  \mu^2 \frac{\kappa(H_1)(z/i)^{1-2H_1}}{(z/i)^2+\beta^2}+\eps \kappa(H_2) (z/i)^{1-2H_2},  \quad \mathrm{Im}(z)>0,
$$ 
where $\kappa(H)=\Gamma(2H+1)\sin (\pi H)$. This definition is extended to the lower half plane through conjugation: 
$$
\Lambda(z;H_1,H_2) := \overline{\Lambda(\overline z;H_1,H_2)}, \qquad \mathrm{Im}(z)<0.
$$
The function in \eqref{spd} is the spectral density of the ``derivative'' of the observation process $Y$ in the stationary case $\beta<0$.
Curiously, it continues to play a role in the non-stationary case as well, albeit in a somewhat implicit way, namely, through its extension to the complex plane.

The so defined extension is sectionally holomorphic on $\mathbb C\setminus \Real$ with finite limits across the cut:
\begin{equation}\label{Llim}
\Lambda^\pm(t;H_1,H_2) := \lim_{z\to t^\pm } \Lambda(z;H_1,H_2), \quad t\in \Real,
\end{equation}
where plus and minus stand for the limits taken in the upper and lower half planes, respectively.
The restriction of $\Lambda(z;H_1,H_2)$ to the imaginary axis coincides with the  spectral density \eqref{spd}:
$$
\Lambda(i\lambda;H_1,H_2)=f(\lambda), \quad \lambda\in \Real.
$$

The configuration of zeros of $\Lambda(z;H_1,H_2)$, which, in view of its definition, must appear in conjugate pairs, plays an important role in the analysis.  
For $H_1 > H_2$, there are two zeros $z_0$ and $\overline z_0$ in the first and the fourth quadrants, respectively. 
As $H_1$ approaches $H_2$, the zeros move  towards positive real semiaxis, and at $H_1=H_2$, collapse to the single 
real zero $t_0 = \sqrt{\beta^2+\mu^2/\eps}$.  For $H_1<H_2$, the zeros disappear and $\Lambda(z;H_1,H_2)$ remains non-vanishing on the cut plane.

The steady-state limit \eqref{limerr} is derived in \cite{AChKM22} as an explicit but rather cumbersome expression in general. 
In a number of meaningful cases, this expression takes a much simpler form. One such case is $H_1=H_2=H$, for which
\eqref{KLeB} is rederived for the full range of $H\in (0,1)$, where the appearance of the real root $t_0$ can be clearly identified.   

Another case corresponds to estimation of the fractional signal in the white noise, that is, when $H_1:=H\in (0,1)$ and $H_2=\frac 1 2$.
Let us define the function, cf. \eqref{Llim},
$$
\theta(t) = \arg\big(\Lambda^+(t; H, \tfrac 1 2)\big), \quad t\in \Real_+,
$$
where the argument branch is chosen so that $\theta(t)$ is 
continuous on $(0,\infty)$. This function has an explicit though somewhat cumbersome formula and it is integrable on $\Real_+$. 
The steady-state limit in this case is found to be
\begin{equation}\label{Pinfty1}
P_\infty\Big(\beta, \frac {\mu}{\sqrt \eps}\Big) = \frac \eps {\mu^2} \left(\frac 1 \pi \int_0^\infty  \theta(t)  dt +\beta  
+
\begin{rcases}   
\begin{dcases}
2\Re (z_0) & \text{if \ } H>\tfrac 1 2 \\
0 & \text{if \ }H<\tfrac 1 2 
\end{dcases}
\end{rcases}
\right),
\end{equation} 
where $z_0$ is the unique zero of $\Lambda(z; H, \frac 1 2)$ in the first quadrant.

In the stationary case $\beta<0$, \eqref{Pinfty1} coincides with the 
well-known expression from the spectral theory of estimation, \cite{Roz67}:
\begin{equation}\label{sfla}
P_\infty\Big(\beta, \frac {\mu}{\sqrt \eps}\Big) = \frac {\eps} {\mu^2} \frac 1 {2\pi} \int_{-\infty}^\infty \log \bigg(1+ \frac{\mu^2}\eps  \frac{\kappa(H)|\lambda|^{1-2H}}{\beta^2+\lambda^2}\bigg) d\lambda. 
\end{equation}
Note that while \eqref{sfla} remains well defined for $\beta\ge 0$ as well, 
it ceases to provide the correct  error limit \eqref{Pinfty1}.

In view of \eqref{epsas}, the exact small noise asymptotics is derived from the steady-state limit \eqref{Pinfty1}:
\begin{equation}\label{Peps0}
P_T\Big(\beta, \frac {\mu}{\sqrt \eps}\Big) =  \frac{\kappa(H)^{1/(2H+1)}}{\sin \frac \pi {2H+1}} (\eps/\mu^2)^{ 2H/(2H+1)}\big(1+o(1)\big), \quad \text as \ \eps\to 0.
\end{equation} 
An alternative way to derive \eqref{Peps0} is by means of the asymptotic approximation for the eigenvalues and eigenfunctions of the covariance operator 
of the fractional OU process, see \cite{ChKM-AiT}.

As mentioned above, such results can be useful in designing simpler filtering algorithms. For example,  
the rate $\eps^{2H/(2H+1)}$ in \eqref{Peps0} coincides with the optimal minimax rate in the nonparametric estimation problem of 
$H$-Holder continuous deterministic signals in small white noise, see \cite{Tsybakov}. Consequently, standard nonparametric procedures, such as 
kernel estimators, etc.,  achieve the optimal rate in the above filtering problem. 

\subsection{Estimation of the Hurst parameter}\label{ss:H}
Estimation of the Hurst parameter of the fBm is an important prototypical problem in statistical inference of time series with 
long range dependence, see \cite{CV10}. In its basic form, the problem consists of constructing estimators for the Hurst exponent $H\in (0,1)$
and the volatility parameter $\sigma\in \Real_+$ from a sample of the scaled fBm $X=\sigma B^H$. 
Reasonable estimators should satisfy performance guarantees, such as consistency in a relevant asymptotic regime. Preferably, they should achieve the best possible rate or, even better, attain the asymptotic risk lower bound.

Both the parameters can be recovered exactly from  
the continuos time sample $X^T = (\sigma B^H_t, t\in [0,T])$ on the horizon of any finite length $T>0$. However, a meaningful statistical model is obtained 
by discretizing the data on a grid of points: 
\begin{equation}\label{dt}
X^{T,\Delta} = \big(\sigma B^H_{\Delta},...,\sigma B^H_{n\Delta}\big),
\end{equation}
where $n = [T/\Delta]$. If both parameters $\sigma$ and $H$ are unknown, they must be estimated {\em jointly}, that is, using only the sample.
If one of the parameters is known, the other can be estimated {\em separately}, that is, its estimator can use the value of the known parameter
as well.

Estimators for  model \eqref{dt} can be studied in the {\em large sample} asymptotic regime with a fixed $\Delta>0$ and $T\to\infty$.
In this case the optimal minimax rate is $T^{-1/2}$, both in the joint or separate estimation problems.
This rate is attainable by the MLE or Whittle's estimator, which are,
in fact, asymptotically efficient in Le Cam's local minimax sense, \cite{D89, D89c}.

An interesting discrepancy occurs in the {\em high frequency}  regime when  $T>0$ is fixed and $\Delta(n):= T/n\to 0$.
In this setting, $H$ and $\sigma$ can be estimated separately at the respective rates 
$$
n^{-1/2} \frac 1{\log n}\quad \text{and}\quad n^{-1/2},
$$
which are minimax optimal, see \cite{K13} and references therein.
These rates are achievable by a number of estimators,  \cite{IL97}, \cite{KW97}, \cite{C01}.
In the joint estimation problem, however, the optimal minimax rates degrade by a logarithmic factor, see \cite{BF17}:
$$
n^{-1/2} \quad \text{and}\quad n^{-1/2}\log n.
$$

A more complicated statistical model is obtained if the sample is contaminated by noise.
This can be done in two different ways: the noise can be added either before or after discretization, see \cite{AC10}.
The first possibility was explored in \cite{GH07}, where the authors considered estimation of $\sigma$ and $H$ from the sample, cf. \eqref{dt}, 
$$
X^{T,\Delta} = \big(\sigma B^H_{\Delta}+\xi_{1},...,\sigma B^H_{n\Delta}+\xi_{n}\big),
$$
where $\xi_1,..., \xi_n$ are  i.i.d. random variables independent of $B^H$. The  rates  
\begin{equation}\label{HGrates}
n^{-1/(4H+2)}\quad \text{and}\quad n^{-1/(4H+2)}\log n
\end{equation}
are shown in \cite{GH07} to be minimax optimal in the joint estimation problem, see also \cite{Sz22}, and 
attainable by adaptive wavelets based estimators. It can be seen that these rates improve by a logarithmic factor if the parameters are  estimated 
separately, as in the noiseless case.

The other natural way to add the noise is to consider the mixed fBm, cf. \eqref{mfBm},
\begin{equation}\label{ctm}
X_t = \sigma B^H_t + \sqrt {\eps} B_t, \quad t\in \Real_+,
\end{equation}
where $B$ is an independent Brownian motion and $\eps>0$ is the known noise intensity parameter. 
The formal derivative of this process is the basic noise model in engineering applications and the estimation problem under consideration 
corresponds to calibration of the parameters of its fractional component, see e.g. \cite{WO92}. As mentioned in Section \ref{sec:mfBm}, 
the measures induced by 
\begin{equation}\label{smpls}
(\sigma B^H_t + \sqrt {\eps} B_t, t\in [0,T])
\quad \text{and} \quad 
(\sqrt \eps B_t, t\in [0,T])
\end{equation}
are mutually absolutely continuous if and only if $H>\frac 3 4$. 
This implies that for such values of $H$ consistent estimation is no longer possible in the high frequency regime
(see, however, \cite{DMS14} for the complementary case $H\le \frac 34$).

The two relevant asymptotics for the continuous time model \eqref{ctm} with $H>\frac 34$ are the {\em large sample} regime with 
a fixed $\eps>0$  and $T\to\infty$ and the {\em small noise} regime with $\eps\to 0$ and a fixed $T>0$.
It was proved in \cite{ChK23} that the optimal minimax rate in the large sample regime remains to be $T^{-1/2}$, both in the joint and separate 
estimation problems. In the small noise regime, the parameters $H$ and $\sigma$ are jointly estimable at the optimal minimax rates, cf.  \eqref{HGrates},
$$
\eps^{1/(4H-2)}\quad \text{and}\quad \eps^{1/(4H-2)}\frac 1{\log \eps^{-1}}.
$$
These rates degrade by a logarithmic factor when the parameters are estimated separately. 

The results in \cite{ChK23} are proved by showing that the model \eqref{ctm} satisfies the local asymptotic normality (LAN) property
with respect to both $T\to\infty$ and $\eps\to 0$. A comprehensive account of the asymptotic theory of estimation can be found in \cite{IKh}.
In a nutshell, in models which satisfy the LAN property, local minimax risks of arbitrary estimators satisfy the Hajek - Le Cam asymptotic 
lower bound, which can be used to identify the optimal minimax estimation rate.
Establishing the LAN property amounts to finding a suitable scaling of the likelihood ratio  under which its logarithm converges to a certain quadratic  
form with a normally distributed coefficient, see \eqref{LAN} below.

The likelihood function for the mixed fBm sample \eqref{ctm} on $[0,T]$ 
with  two dimensional parameter $\theta := (H,\sigma)\in (\frac 3 4,1)\times (0,\infty)$ 
is the Radon-Nikodym density between measures induced by the processes in \eqref{smpls}. It can be derived from the canonical 
innovation representation in \cite{Hi68}:
$$
L_T(\theta):= \frac{d\P^T_\theta}{d\P^T}(X^T) =\exp\left(\frac 1 \eps \int_0^T \rho_t(X,\theta)dX_t
-\frac 1{2\eps}   \int_0^T \rho_t(X,\theta)^2dt\right),
$$
where 
$$
\rho_t(X,\theta) = \int_0^t g(t,t-s;\theta)dX_s.
$$
The function $g(t,s;\theta)$ solves the integral equation of the second kind, cf. \eqref{kind2},
\begin{equation}\label{eq}
\eps g(t,s;\theta) + \int_0^t  K_{\theta}(r-s) g(t,r;\theta) dr =   K_{\theta}(s), \quad 0<s<t,
\end{equation}
with the kernel 
$$
K_\theta (\tau)= \sigma^2 H(2H-1) |\tau|^{2H-2}, \quad \tau \in \Real.
$$

The LAN property with respect to $T\to \infty$ requires finding an appropriate scaling function $\phi(T, \theta_0)$ for which the likelihood ratio  satisfies the 
decomposition 
\begin{equation}\label{LAN}
\log \frac{L_T\big(\theta_0 + \phi(T, \theta_0) u\big)}{L_T(\theta_0 )} = u^\top Z_{T} -\frac 1 2 \|u\|^2 + r_T, \quad \forall u\in \Real^2,
\end{equation}
such that  
$$
Z_T \xrightarrow[T\to\infty]{\mathcal{L}(P_{\theta_0})}  N(0,\mathrm{Id})\quad \text{and}\quad   r_T \xrightarrow[T\to\infty]{ P_{\theta_0} } 0.
$$ 
The likelihood ratio is a particular functional of the solution to \eqref{eq}. Thus establishing the LAN property \eqref{LAN}
again reduces to asymptotic analysis of an integral equation.

The appropriate choice in the large sample regime $T\to\infty$ is  $\phi(T, \theta_0)= \frac 1 {\sqrt T} I(\theta_0,\eps)^{-1/2}$  
where 
\begin{equation}\label{Ieps}
I(\theta,\eps)  = \frac 1 {4\pi} \int_{-\infty}^\infty \nabla^\top \log \big(\eps +  \widehat K_\theta(\lambda)\big)
\nabla \log \big(\eps +  \widehat K_\theta(\lambda)\big) d\lambda>0
\end{equation}
is the Fisher information rate matrix ($\nabla$ stands for the gradient with respect $\theta$). This is, again, 
the analog of  Whittle's formula \cite{W53, W62} in continuous time. 

The LAN property with respect to $\eps\to 0$ is shown to hold for any scaling function $\phi(\eps,\theta_0)$ 
which satisfies 
\begin{equation}\label{condphi}
\phi(\eps,\theta_0)^\top M(\eps,\theta_0) T I(\theta_0,1)M(\eps,\theta_0)^\top\phi(\eps,\theta_0)\xrightarrow[\eps\to 0]{}\mathrm{Id},
\end{equation}
with $I(\theta_0,1)$ defined in \eqref{Ieps}
and  
$$
M(\eps, \theta) = \eps^{-1/(4H-2)}\begin{pmatrix*}[c]
1 \ \ & -2\sigma^2  \log\eps^{-1/(2H-1)}\\
0 \ \ & 1
\end{pmatrix*}.
$$
Condition \eqref{condphi} cannot be satisfied by any diagonal matrix $\phi(\eps,\theta_0)$, since in this case the limit, if it exists and finite, 
must be a singular matrix. Otherwise the choice of $\phi(\eps,\theta_0)$ is not unique.   
The upper and lower triangular Cholesky factors of the matrix $M(\eps,\theta_0) I(\theta_0,1)M(\eps,\theta_0)^\top$ yield the optimal 
minimax estimation rates for $H$ and $\sigma^2$ in the joint estimation problem. The optimal rates in the separate estimation are found by 
considering the restrictions of this LAN property to the one-dimensional families corresponding to each one of the parameters.

\section{About the method}\label{sec:3}
In this section we outline the approach to asymptotic analysis of integral equations, such as those encountered in the 
previous sections. It is inspired by the technique pioneered in \cite{Ukai} and \cite{P74,P03}, which was originally 
developed for asymptotic spectral approximation for operators with weakly singular difference kernels. 

\subsection{The roadmap}

In a nutshell, the main idea is to reduce the problem to a certain system of integro-algebraic equations which depend on a large parameter. 
This system turns out to be more amenable to asymptotic analysis as the parameter tends to infinity. 
Though the implementation of the method is somewhat specific to concrete problems, the main steps do follow a common pattern.

\begin{enumerate}

\item The Laplace transform of the solution $\psi(\cdot)$  to an equation at  hand 
\begin{equation}\label{Lap}
\widehat \psi(z) =  \int_{0}^1 \psi(x) e^{-zx }dx, \quad z\in \mathbb C
\end{equation}
is shown to have a representation of the form
\begin{equation}\label{rep}
\widehat \psi(z) = \frac{\Phi_0(z) + e^{-z}\Phi_1(-z)}{\Lambda(z)}, \quad z\in \mathbb C.
\end{equation}
Here $\Lambda(z)$ is an explicit function, sectionally holomorphic on the cut plane $\mathbb C\setminus \Real$, that is, holomorphic 
everywhere except for the real line, where it has a jump discontinuity.
The functions $\Phi_0(z)$ and $\Phi_1(z)$ are defined by certain functionals of the solution to the equation being considered and 
thus are treated as unknowns. They are sectionally holomorphic on $\mathbb C\setminus \Real_+$ and their 
growth at zero and at infinity must be known. Usually, these functions grow at infinity as polynomials of a finite degree 
with possibly unknown coefficients. 

\item Analyticity of the Laplace transform \eqref{Lap} on the whole complex  plane implies continuity across $\Real$ of the expression in the right 
hand side of \eqref{rep}:
\begin{equation}\label{limlim}
\lim_{z\to t^+}\frac{\Phi_0(z) + e^{-z}\Phi_1(-z)}{\Lambda(z)}
=
\lim_{z\to t^-}\frac{\Phi_0(z) + e^{-z}\Phi_1(-z)}{\Lambda(z)}, \quad t\in \Real.
\end{equation}
Here plus and minus under the limits mean that $z$ approaches $t$ in the upper and lower half planes respectively.
Analyticity of  \eqref{Lap} also implies the algebraic conditions
\begin{equation}\label{poles}
\Phi_0( z_j) + e^{- z_j}\Phi_1(- z_j) =0, \quad j=0,...,m
\end{equation}
where $z_j$'s are  zeros of $\Lambda(z)$.

\item  Using the technique for solving Hilbert's boundary value problems, see \cite{Gahov}, condition \eqref{limlim} can be reduced to 
integral equations of the form
\begin{equation}\label{SDeq}
\begin{aligned}
S(t) = & \ \phantom+\ \frac 1 \pi \int_0^\infty \frac{h_0(\tau)}{\tau +t}e^{-\nu \tau}S(\tau)d\tau + P_S(t) 
\\
D(t) = & -\frac 1 \pi \int_0^\infty \frac{h_0(\tau)}{\tau +t}e^{-\nu \tau}D(\tau)d\tau +  P_D(t)
\end{aligned}\quad t\in \Real_+,
\end{equation}
where $h_0(\cdot)$ is a known function and $\nu$ is a large parameter. 
The solutions $S(\cdot)$ and $D(\cdot)$ are related to $\Phi_0(\cdot)$ and $\Phi_1(\cdot)$ by a specific linear transformation.
$P_S(\cdot)$ and $P_D(\cdot)$ are polynomials whose degrees are chosen to match the growth of  $S(\cdot)$ and $D(\cdot)$ at infinity.

\item The coefficients of  polynomials $P_S(\cdot)$ and $P_D(\cdot)$ are determined by the conditions which guarantee nontrivial solution to 
the integro-algebraic system \eqref{SDeq} and \eqref{poles}. 

\item Asymptotic analysis is based on convergence of the integral terms in \eqref{SDeq} to zero as $\nu\to\infty$.

\end{enumerate}

\subsection{Some details}

Let us elaborate on the above program using, as an example, eigenvalue 
problem \eqref{eig} for the fBm with covariance \eqref{cov}: 
\begin{equation}\label{thepr}
\int_0^1 \frac 1 2 \Big(|t|^{2-\alpha} +   |s|^{2-\alpha} - |t-s|^{2-\alpha} \Big)\phi(s)ds = \lambda \phi(t), \quad t\in [0,1],
\end{equation} 
where we defined $\alpha := 2-2H$ for convenience. We will consider the case $H>\frac 1 2$, that is, $\alpha\in (0,1)$. 

\subsubsection{Reduction to a difference kernel}
The starting point is reduction to an equivalent problem which involves a difference kernel. For equation \eqref{thepr}, this can be done by switching to 
the function $\psi(x) = \int_x^1 \phi(t)dt$. A direct calculation shows that $\psi(\cdot)$ solves the generalized eigenvalue problem
\begin{equation}
\label{geig}
\begin{aligned}
\int_0^1  & c_\alpha  |x-y|^{-\alpha}  \psi(y)dy = -\lambda \psi''(x), \quad x\in [0,1] \\
&
\psi(1)=0, \; \psi'(0)=0
\end{aligned}
\end{equation}
where $c_\alpha : = (1- \alpha/2)(1-\alpha)$. Such reduction is case specific, see, e.g., \cite{ChK21}, \cite{ChKM-AiT}.

\subsubsection{Representation of the Laplace transform}
The next step is to derive  representation \eqref{rep}.
For equation \eqref{geig}  it holds with 
\begin{equation}\label{Lambda}
\Lambda(z) = \frac 1 z \frac{\Gamma(\alpha)}{c_\alpha} \Big(\lambda z^2 + Q(z)\Big), \quad z\in \mathbb C\setminus \Real,
\end{equation}
where $Q(\cdot)$ is a sectionally holomorphic extension to $\mathbb C\setminus \Real$ of the Fourier transform of the kernel in \eqref{geig}: 
$$
\chi(\omega) = \int_{-\infty}^\infty c_\alpha |\tau|^{-\alpha} e^{-i\omega \tau}d\tau = d_\alpha |\omega|^{\alpha-1}, \quad \omega \in \Real\setminus\{0\},
$$
where we defined the constant
$$
d_\alpha = \dfrac{  c_\alpha }{\Gamma(\alpha)  }\dfrac{\pi   }{ \cos (\frac \pi2 \alpha)}.
$$
Such extension is given by the closed form formula
$$
Q(z) =  \begin{cases}
d_\alpha(z/i)^{\alpha-1}, & \Im(z)>0, \\
d_\alpha(-z/i)^{\alpha-1}, & \Im(z)<0.
\end{cases}
$$
Clearly,  $Q(i\omega)=\chi(\omega)$ for $\omega \in \Real$ and $Q(z)$ has a jump discontinuity on the real line.    
The applicability of the method is primarily determined by  existence of such an extension.

\subsubsection{The equivalent Hilbert problem}

In view of the symmetries of $\Lambda(z)$, condition \eqref{limlim} reduces to 
\begin{equation}
\label{Hp}
\begin{aligned}
&
\Phi_0^+(t) - e^{2i\theta(t)}\Phi_0^-(t) = 2i  e^{-t} e^{i\theta(t)}\sin\theta(t) \Phi_1(-t)
\\
&
\Phi_1^+(t) - e^{2i\theta(t)}\Phi_1^-(t) = 2i e^{-t}  e^{i\theta(t)}\sin\theta(t) \Phi_0(-t)
\end{aligned}\quad t\in \Real_+,
\end{equation}
where we defined $\theta(t):=\arg(\Lambda^+(t))$ and denoted   
$
\Phi_j^\pm(t) = \lim_{z\to t^\pm}\Phi_j(z).
$
The function $\Lambda(z)$ in \eqref{Lambda} has a pair of purely imaginary zeros $\pm z_0 = \pm  i \nu$ where  
$\nu$ is the large parameter, related to $\lambda$ through the formula 
\begin{equation}\label{nuflan}
\nu^{\alpha-3} =  \frac{\lambda}{d_\alpha}.
\end{equation}
Condition \eqref{poles} reduces to  
\begin{equation}\label{z0}
\Phi_0( z_0) + e^{- z_0}\Phi_1(- z_0) =0.
\end{equation}

Equations \eqref{Hp} can be viewed as boundary conditions for a pair of unknown functions $\Phi_0(z)$ and $\Phi_1(z)$, sectionally holomorphic 
on $\mathbb C\setminus \Real_+$. Thus we arrive at a variation of the Hilbert boundary value problem of finding all such functions, which, in addition, 
must satisfy the a priori growth conditions and the algebraic constraint \eqref{z0}. Note that $\Phi_0(z)\equiv 0$ and 
$\Phi_1(z)\equiv 0$ are solutions, which correspond to the trivial solutions of problem \eqref{thepr}. As we will see below, 
the nontrivial solutions exist only for special values of $\nu$. The corresponding values of $\lambda$, given by the formula \eqref{nuflan},
are precisely the eigenvalues and,  for each such value of $\lambda$, the corresponding eigenfunction can be reconstructed by plugging $\Phi_0(z)$ and $\Phi_1(z)$ into \eqref{rep} and inverting the Laplace transform.

\subsubsection{Integro-algebraic system}
The Hilbert problem for  boundary conditions \eqref{Hp} can be reduced to a system \eqref{SDeq}. Solutions to this system are 
linear combinations with coefficients $a_j$ and $b_j$ of solutions to the auxiliary stand alone equations:
\begin{equation}
\label{pq}
\begin{aligned}
& p_\pm(t  ) = \pm \frac 1 {\pi } \int_0^\infty \frac{h_0(s)e^{-\nu s} }{s +t }p_\pm (s) ds  + 1, \\
& q_{\pm }(t  ) = \pm \frac 1 {\pi } \int_0^\infty \frac{h_0(s)e^{-\nu s} }{s+t  }q_\pm (s) ds + t,
\end{aligned}\quad\quad t\in \Real_+,
\end{equation}
where $h_0(\cdot)$ is an explicit function and $\nu$ defined in \eqref{nuflan}.
These equations are derived using standard tools from complex analysis, such as the Sokhotski–Plemelj theorem, see \cite{Gahov}.
It can be shown that, at least for all sufficiently large $\nu$, the integral operators in \eqref{pq} are contractions on $L_2(\Real_+)$  
and therefore unique solutions exist.

Let us extend the domain of $p_\pm(\cdot)$ and $q_\pm(\cdot)$  
by replacing $t\in \Real_+$ with $z\in  \mathbb C\setminus \Real_{-}$ in \eqref{pq} and define 
\begin{equation}\label{abpm}
\begin{aligned}
& a_\pm(z) = p_+(z)\pm p_-(z), \\
& b_\pm(z) = q_+(z)\pm q_-(z).
\end{aligned}
\end{equation}
It can be shown that the functions, which satisfy the boundary conditions \eqref{Hp} and the relevant a priori growth estimates, are given by
\begin{equation}
\label{PhPhA}
\begin{aligned}
& \Phi_0(z\nu) =  
c_2 \nu X_0(z)\big(b_+(-z)-b_\alpha a_+(-z) \big)+c_1 X_0(z) a_-(-z),
\\
& \Phi_1(z\nu) = 
c_2 \nu  X_0(z) \big(b_-(-z)-b_\alpha a_-(-z)\big)+c_1 X_0(z) a_+(-z),
\end{aligned}
\end{equation}
where the constant $b_\alpha$ and the function $X_0(\cdot)$ are explicit and $c_1$ and $c_2$ are unknown parameters (which determine the coefficients 
of $P_S(\cdot)$ and $P_D(\cdot)$ in \eqref{SDeq}).  

If we now plug the expressions \eqref{PhPhA} into condition \eqref{z0}, we arrive at a homogeneous two-dimensional system of linear 
algebraic equations 
\begin{equation}\label{sys}
\eta c_1+ \nu \xi c_2  =0
\end{equation}
for the real unknowns $c_1$ and $c_2$
with the complex valued coefficients 
\begin{equation}\label{xieta}
\begin{aligned}
&
\xi := e^{ i\nu/2} X_0(i)\big(b_+(-i)-b_\alpha a_+(-i) \big)+e^{-i\nu/2} X_0(-i) \big(b_-(i)-b_\alpha a_-(i)\big), \\
&
\eta:=  e^{i\nu/2} X_0(i) a_-(-i)+e^{-i\nu/2} X_0(-i) a_+(i).
\end{aligned}
\end{equation}

The trivial solution $c_1=c_2=0$  to \eqref{sys} correspond to the trivial solution  to eigenvalue problem \eqref{thepr}. 
Nontrivial solutions exist if and only if the determinant of the coefficients matrix vanishes. Since $c_1$ and $c_2$ are real, this is equivalent to 
the condition
\begin{equation}\label{alg}
\Im\{\xi \overline {\eta}\} = 0.
\end{equation}
The quantities $\xi$ and $\eta$, defined in \eqref{xieta}, are determined by the solutions to \eqref{pq}, which in turn, depend on $\nu$. 
Thus, we are faced with an integro-algebraic system of equations which consists of \eqref{pq} and \eqref{alg}. 

Let $(p_\pm, q_\pm, \nu)$ be a solution to this system. Plugging $\nu$ into \eqref{nuflan} immediately gives an eigenvalue $\lambda$ to 
\eqref{thepr}. Moreover, for this value of $\nu$, equations \eqref{sys} give the linear relation $c_1 = \nu \xi /\eta c_2$. 
Plugging it into \eqref{PhPhA} along with the corresponding solutions $p_\pm(\cdot)$ and $q_\pm(\cdot)$ to \eqref{pq} and \eqref{abpm}  
we obtain $\Phi_0(\cdot)$ and $\Phi_1(\cdot)$,
where the only unknown is the multiplicative factor $c_2$. 
Finally, substituting these functions into \eqref{rep} and inverting the Laplace transform we get $\psi(\cdot)$ and consequently the 
eigenfunction $\phi(\cdot)=-\psi'(\cdot)$. The factor $c_2$ can be eliminated, e.g., by normalizing to unit norm.

\subsubsection{Asymptotic analysis}

The equivalent integro-algebraic system \eqref{pq} and \eqref{alg} does not appear, at the outset, to be any simpler than 
the original eigenproblem \eqref{thepr} and there is little hope for reasonably explicit solutions. 
Remarkably, it turns out to be more tractable as $\nu\to\infty$. It can be argued to have countably many solutions, 
and all but a finite number of them can be approximated asymptotically. 
The integral term in \eqref{pq} is asymptotically negligible as $\nu\to \infty$.
More specifically, for any fixed $z\not \in \Real_-$,
\begin{align*}
|p_{\pm}(z) - 1| \le C\nu^{-1}   \quad \text{and}\quad |q_{\pm}(z) -z | \le C \nu^{-2},
\end{align*}
where constant $C$ may depend on $z$. Consequently, cf. \eqref{abpm},
\begin{align*}
a_+(\pm i) =\, &  2 + O(\nu^{-1}), \\
b_+(\pm i) =& \pm 2i + + O(\nu^{-2}), \\
a_-(\pm i) =&    O(\nu^{-1}), \\
b_-(\pm i) =&   O(\nu^{-2}), \\
\end{align*}
and, cf. \eqref{xieta},
$$
\begin{aligned}
&
 \xi  =   e^{ i\nu/2} X_0(i)\big(-2i-2b_\alpha   \big) + O(\nu^{-1}),  \\
&
\eta =   2 e^{-i\nu/2} X_0(-i)  + O(\nu^{-1}).
\end{aligned}
$$
Hence, asymptotically,  \eqref{alg}  reduces to 
$$
\mathrm{Im} \Big(e^{ i\nu } X_0(i)\big( i+b_\alpha   \big)   \overline{X_0(-i)}\Big) =O(\nu^{-1}), \quad \nu\to \infty,
$$
and, therefore,
\begin{equation}\label{nunu}
\nu_n + \pi n +  \arg{X_0(i)} - \arg{X_0(-i)} + \arg(i+b_\alpha) + O(n^{-1}), \quad n\to\infty,
\end{equation}
where all terms can be found in closed form. This and \eqref{nuflan} give \eqref{nun}.

Any  $\nu_n$ in \eqref{nunu}, for sufficiently large $n$, corresponds to an eigenvalue in  problem \eqref{thepr}.
Thus we found all but finitely many entries in the ordered sequence of eigenvalues. However, the enumeration in \eqref{nunu} may be shifted by a finite 
integer with respect to the enumeration which puts {\em all} eigenvalues in the decreasing order. 
This shift can be identified by a calibration procedure based on the exact formula 
in the Brownian case $\alpha =1$. The corresponding approximations \eqref{fBmeigfn} for the eigenfunctions are obtained by asymptotic analysis of the inversion formulas 
of the Laplace transform.

%\bibliographystyle{spmpsci}
%\bibliography{/Users/Pavel/Dropbox/Pasha_Masha/bibliography/fBm}

\def\cprime{$'$} \def\cprime{$'$} \def\cydot{\leavevmode\raise.4ex\hbox{.}}
  \def\cprime{$'$} \def\cprime{$'$} \def\cprime{$'$}

\end{document}